\documentclass[graybox]{svmult}

\usepackage{mathptmx}       % selects Times Roman as basic font
\usepackage{helvet}         % selects Helvetica as sans-serif font
\usepackage{courier}        % selects Courier as typewriter font
\usepackage{type1cm}        % activate if the above 3 fonts are
\usepackage{makeidx}         % allows index generation
\usepackage{graphicx}        % standard LaTeX graphics tool
\usepackage{multicol}        % used for the two-column index
\usepackage[bottom]{footmisc}% places footnotes at page bottom

\makeindex             % used for the subject index
\begin{document}
\title*{Local Isometric immersions of pseudo-spherical surfaces and evolution equations}
% Use \titlerunning{Short Title} for an abbreviated version of
% your contribution title if the original one is too long
\author{\author{Nabil Kahouadji, Niky Kamran and Keti Tenenblat}}
% Use \authorrunning{Short Title} for an abbreviated version of
% your contribution title if the original one is too long
\institute{Nabil Kahouadji \at Department of Mathematics, Northwestern University, \email{nabil@math.northwestern.edu}
\and Niky Kamran\at Department of Mathematics and Statistics, McGill University, \email{nkamran@math.mcgill.ca}
\and Keti Tenenblat\at Departamento de Matem\`atica, Universidade de Brasilia, \email{keti@mat.unb.br}}
%
% Use the package "url.sty" to avoid
% problems with special characters
% used in your e-mail or web address
%
\maketitle
\centerline{\it To Walter Craig, with friendship and admiration}
\abstract*{Each chapter should be preceded by an abstract (10--15 lines long) that summarizes the content. The abstract will appear \textit{online} at \url{www.SpringerLink.com} and be available with unrestricted access. This allows unregistered users to read the abstract as a teaser for the complete chapter. As a general rule the abstracts will not appear in the printed version of your book unless it is the style of your particular book or that of the series to which your book belongs.
Please use the 'starred' version of the new Springer \texttt{abstract} command for typesetting the text of the online abstracts (cf. source file of this chapter template \texttt{abstract}) and include them with the source files of your manuscript. Use the plain \texttt{abstract} command if the abstract is also to appear in the printed version of the book.}

\abstract{The class of differential equations describing pseudo-spherical surfaces, first introduced by Chern and Tenenblat~\cite{CT}, is characterized by the property that to each solution of a differential equation, within the class, there corresponds a 2-dimensional Riemannian metric of curvature equal to $-1$. The class of differential equations describing pseudo-spherical surfaces carries close ties to the property of complete integrability, as manifested by the existence of infinite hierarchies of conservation laws and associated linear problems. As such, it contains many important known examples of integrable equations, like the sine-Gordon, Liouville and KdV equations. It  also gives rise to many new families of integrable equations. The question we address in this paper concerns the local isometric immersion of pseudo-spherical surfaces in ${\bf E}^{3}$ from the perspective of the differential equations that give rise to the metrics. Indeed, a classical theorem in the differential geometry of surfaces states that any pseudo-spherical surface can be locally isometrically immersed in ${\bf E}^{3}$. In the case of the sine-Gordon equation, one can derive an expression for the second fundamental form of the immersion that depends only on a jet of finite order of the solution of the pde. A natural question is to know if this remarkable property extends to equations other than the sine-Gordon equation within the class of differential equations describing pseudo-spherical surfaces. In an earlier paper~\cite{KKT}, we have shown that this property fails to hold for all other second order equations, except for those belonging to a very special class of evolution equations. In the present paper, we consider a class of evolution equations 
for $u(x,t)$ of order $k\geq 3$ describing pseudo-spherical surfaces. We show that whenever an isometric immersion in ${\bf E}^3$ exists, depending on a jet 
of finite order of $u$, then the coefficients of the second fundamental forms 
are functions of the independent variables $x$ and $t$ only.   }
\section{Introduction and Statement of Results}
\label{sec:1}

The notion of a partial differential equation describing pseudo-spherical surfaces was defined and studied extensively in a paper by Chern and Tenenblat~\cite{CT}. The class of these equations is of particular interest because it enjoys a remarkable set of integrability properties in the case when a parameter playing the role of a spectral parameter is present in the $1$-forms associated to the pseudo-spherical structure. Indeed, one obtains in that case an infinite sequence of conservation laws and  an associated linear problem whose integrability condition is the given partial differential equation. \footnote{It is worth pointing out at this stage that the conservation laws arising from the geometry of pseudo-spherical surfaces may be non-local. We refer to \cite{R106} and the references therein for an explicit treatment of the relationship between the conservation laws obtained in \cite{CaT} and the standard series of conservation laws obtained via the classical Wahlquist-Estabrook construction. We also refer to \cite{H},\cite{Sa06} and \cite{Sa07} for the study of the integrability properties of some specific families of equations describing pseudo-spherical surfaces.} We begin by recalling some basic definitions.  A partial differential equation 
\begin{equation}\label{pde}
\Delta(t,x,u,u_{t},u_{x},\ldots,u_{t^{l}x^{k-l}})=0,
\end{equation}
is said to describe pseudo-spherical surfaces if there exist $1$-forms 
\begin{equation}\label{forms}
\omega^{i}=f_{i1}dx+f_{i2}dt,\quad 1\leq i \leq3,
\end{equation}
where the coefficients $f_{ij},\,1\leq i \leq 3,\,1\leq j\leq 2,$ are smooth functions of $t,x,u$ and finitely many derivatives of $u$, such that the structure equations for a surface of Gaussian curvature equal to $-1$,
\begin{equation}\label{struct}
d\omega^{1}=\omega^{3} \wedge\omega^{2},\quad d\omega^{2}=\omega^{1} \wedge\omega^{3},\quad d\omega^{3}=\omega^{1} \wedge\omega^{2}
\end{equation}
are satisfied if and only if $u$ is a solution of (\ref{pde}) for which $\omega^{1}\wedge \omega^{2}\neq 0$. In other words, for every smooth solution of (\ref{pde}) such that $\omega^{1}$ and $\omega^{2}$ are linearly independent, we obtain a Riemannian metric 
\begin{equation}\label{metric}
ds^{2}=(\omega^{1})^{2}+(\omega^{2})^{2},
\end{equation}
of constant Gaussian curvature equal to $-1$, with $\omega^{3}$ being the Levi-Civita connection $1$-form. This condition is equivalent to the integrability condition for the linear problem given by
\begin{equation}\label{linear}
dv^{1}=\frac{1}{2}(\omega^2{}v^{1}+(\omega^{1}-\omega^{3})v^{2}),\quad dv^{2}=\frac{1}{2} ((\omega^{1}+\omega^{3})v^{1}-\omega^{2}v^{2}).
\end{equation}
For the purposes of this paper, the motivating example of a partial differential equation describing pseudo-spherical surfaces is the sine-Gordon equation
\begin{equation}\label{sG}
u_{tx}=\sin u, 
\end{equation}
for which a choice of $1$-forms (\ref{forms}) satisfying the structure equations (\ref{struct}) is given by
\begin{equation}\label{sGcof1}
\omega^1 = \frac{1}{\eta}\sin u \, dt, \quad \omega^2 = \eta\, dx+\frac{1}{\eta}\cos u \,dt,\quad \omega^3 = u_{x}\,dx,
\end{equation}
where $\eta$ is a non-vanishing real parameter. This continuous parameter is closely related to the parameter appearing in the classical B\"acklund transformation for the sine-Gordon equation and is central to the existence of infinitely many conservation laws for the sine-Gordon equation. It is noteworthy that there may be different choices of $1$-forms satisfying the structure equations (\ref{struct}) for a given differential equation. For example, for the sine-Gordon equation (\ref{sG}), a choice different from the one given in (\ref{sGcof1}) is given by
\begin{equation}\label{sGcof2}
\omega^1 = \cos \frac{u}{2}( dx+dt),\quad \omega^2 =\sin \frac{u}{2} (dx - dt),\quad \omega^3 = \frac{u_{x}}{2} dx - \frac{u_{t}}{2} dt.
 \end{equation}
Partial differential equations (\ref{pde}) which describe pseudo-spherical surfaces and for which one of the components $f_{ij}$ can be chosen to be a continuous parameter will be said to describe $\eta$ pseudo-spherical surfaces. One important feature of the differential equations describing $\eta$-pseudo-spherical surfaces is that each such differential equation is the integrability condition of a linear system of the form (\ref{linear}), which may be used as a starting point in the inverse scattering method and lead to solutions of the differential equation (see for example \cite{BRT}). 

% These equations, which have been completely classified by Chern and Tenenblat in \cite{ChernTenenblat}, %provide an extensive class of integrable non-linear partial differential equations in two independent % variables.

It is therefore an interesting problem to characterize  the class of differential equations describing $\eta$-pseudo-spherical surfaces, and this is precisely what Chern and Tenenblat \cite{CT} did for $k$-th order evolution equations
\begin{equation}\label{evoleq}
u_t= F(u,u_x,..., u_{x^k}).
\end{equation} 
These results were extended to more general classes of differential equations in \cite{R88}, \cite{R89}, \cite{RT90}, and \cite{RT92}.  One can also remove the assumption that $f_{21}=\eta$ and perform a complete characterization  of evolution equations of the form (\ref{evoleq}) 
which describe pseudo-spherical surfaces, as opposed to $\eta$ pseudo-spherical surfaces \cite{KT}. It is also noteworthy that the classification results obtained by Chern and Tenenblat \cite{CT} for $\eta$-pseudo-spherical surfaces were extended in \cite{R98} to differential equations of the form 
$u_t= F(x,t, u,u_x,..., \partial u/\partial x^k)$.  
Finally we mention that the concept of a differential equation that describes pseudo-spherical surfaces has a spherical counterpart \cite{DT}, where similar classification results have been obtained.  Further developments can be found in 
\cite{CaT}, \cite{FT}, \cite{Go}, \cite{GR}, \cite{JT}, \cite{FOR}, \cite{R98},         \cite{R106}.

The property of a surface being pseudo-spherical is by definition intrinsic since it only depends on its first fundamental form. It is only recently \cite{KKT} that the problem has been considered of locally isometrically immersing in ${\bf E}^{3}$ the pseudo-spherical surfaces arising from the solutions of partial differential equations describing pseudo-spherical surfaces. Let us first recall that any pseudo-spherical surface can be locally isometrically immersed into three-dimensional Euclidean space ${\bf E}^{3}$. This means that to any solution $u$ of a partial differential equation (\ref{pde}) describing pseudo-spherical surfaces (for which $\omega^{1}\wedge \omega^{2}\neq 0$), there corresponds a local isometric immersion into ${\bf E}^{3}$ for the corresponding metric of constant Gaussian curvature equal to $-1$. The problem investigated in \cite{KKT} was to determine to what extent the {\em second fundamental form} of the immersion could be expressed in terms of the solution $u$ of a second-order equation and {\emph{finitely many}} of its derivatives. The motivation for this question came from a remarkable property of the sine-Gordon equation, which we now explain.  Let us first derive a set of necessary and sufficient conditions that the components $a,b,c$ of the second fundamental form of any local isometric immersion into ${\bf E}^{3}$ of a metric of constant curvature equal to $-1$ must satisfy. Recall that $a,b,c$ are defined by the relations
\begin{equation}\label{w13_w23}
\omega^3_{1} = a\omega^1+b\omega^2, \quad \omega^3_{2} = b\omega^1+c\omega^2,
\end{equation}
where the $1$-forms $\omega^3_{1}, \omega^3_{2}$ satisfy the structure equations
\begin{equation}\label{Codazzi}
d\omega^3_{1} = -\omega^3_{2}\wedge\omega^2_{1}, \quad d\omega^3_{2} = -\omega^3_{1}\wedge\omega^1_{2},
\end{equation}
equivalent to the Codazzi equations, and the Gauss equation 
\begin{equation}%\label{Gauss}
ac-b^2=-1.
\end{equation}
For the sine-Gordon equation, with the choice of $1$-forms $\omega^{1},\,\omega^{2}$ and $\omega^{3}=\omega^{2}_{1}$ given by (\ref{sGcof2}), it is easily verified that the $1$-forms $\omega^{3}_{1},\omega^{3}_{2}$ are given by
\begin{eqnarray*}
\omega^3_{1} &=& \sin\frac{u}{2} (dx+dt) =  \tan \frac{u}{2}\omega^1,\\
\omega^3_{2} &=& -\cos\frac{u}{2} (dx - dt) = -\cot \frac{u}{2}\omega^2.
\end{eqnarray*}
It is a most remarkable fact that the components $a,b,c$ of the second fundamental form that we have just obtained depend only on the solution $u$ of the sine-Gordon equation. Our main goal in \cite{KKT} was to investigate to what extent this property was true for all second-order equations describing pseudo-spherical surfaces, in the sense that $a,b,c$ should only depend on $u$ and at most {\emph{finitely many derivatives}} of $u$. We showed that this is an extremely rare event, essentially confined to the sine-Gordon equation. Indeed, what we proved was that except for the equation
\begin{equation}
u_{xt}=F(u), \quad F^{''}(u)+\alpha u=0,
\end{equation}
where $\alpha$ is a positive constant, every second-order partial differential equation describing $\eta$-pseudo-spherical surfaces is such that either $a,b,c$ are universal functions of $t,x$, meaning that they are independent of $u$, or they do not factor through any finite-order jet of $u$. The starting point of the proof of this result is a set of necessary and sufficient conditions, given in terms of the coefficients $f_{ij}$ of the $1$-forms (\ref{forms}), for $a,b$ and $c$ to be the components of the second fundamental form of a local isometric immersion corresponding to a solution of (\ref{pde}). These are equivalent to the Gauss and Codazzi equations, and are easily derived.  We consider the pair of vector fields $(e_{1}, e_{2})$ dual to the coframe $(\omega^1, \omega^2)$. It is given by 
\begin{equation}
 \left|\begin{array}{cc}f_{11} & f_{21} \\f_{12} & f_{22}\end{array}\right|e_{1} = f_{22}\partial_{x} - f_{21}\partial_{t}, \quad  \left|\begin{array}{cc}f_{11} & f_{21} \\f_{12} & f_{22}\end{array}\right|e_{2} = -f_{12}\partial_{x} + f_{11}\partial_{t}.
 \end{equation}  
By feeding these expressions into the structure equations (\ref{Codazzi}), we obtain
\begin{eqnarray}
d\omega^3_{1}= \Big( db(e_{1}) -  da(e_{2})\Big)\omega^1\wedge\omega^2 - a\,\omega^2\wedge\omega^3 +b\,\omega^1\wedge\omega^3, \label{dw13}\\ 
d\omega^3_{2}= \Big( dc(e_{1}) -  db(e_{2})\Big)\omega^1\wedge\omega^2 - b\,\omega^2\wedge\omega^3 +c\,\omega^1\wedge\omega^3, \label{dw23}
\end{eqnarray}
Denoting by $D_{t}$ and $D_{x}$ the total derivative operators, these are equivalent to 
\begin{eqnarray}\label{Eq1}
 f_{11} D_{t}a + f_{21}D_{t}b - f_{12}D_{x}a - f_{22}D_{x}b - 2b \left|\begin{array}{cc}f_{11} & f_{31} \\f_{12} & f_{32}\end{array}\right|+ (a-c)\left|\begin{array}{cc}f_{21} & f_{31}\\f_{22} & f_{32}\end{array}\right|   = 0,\\\label{Eq2}
 f_{11} D_{t}b + f_{21}D_{t}c - f_{12}D_{x}b - f_{22}D_{x}c +(a-c) \left|\begin{array}{cc}f_{11} & f_{31} \\f_{12} & f_{32}\end{array}\right|+ 2b\left|\begin{array}{cc}f_{21} & f_{31}\\f_{22} & f_{32}\end{array}\right|  = 0.
 \end{eqnarray}
These differential constraints, which amount to the Codazzi equations, have to be augmented by the Gauss equation 
\begin{equation}\label{Gauss}
ac-b^2=-1.
\end{equation}  
The proof of the main result of \cite{KKT} consists in a detailed case-by-case analysis of the equations (\ref{Eq1}), (\ref{Eq2}) and (\ref{Gauss}), where use is made of the expressions and constraints on the $f_{ij}$'s that result from the classification results of \cite{CT}, \cite{RT90},  and where we assume that $a, b$ and $c$ depend on $t,x,u$ and only finitely many derivatives of $u$. 

Our goal in the present paper is to extend the results of \cite{KKT} concerning the components $a,b,c$ of the second fundamental form to the case of $k$-th order evolution equations. Our main result is given by:

\begin{theorem}\label{LemCoeff}Let 
\begin{equation}\label{evolk}
u_t=F(u, u_{x}, \dots, u_{x^{k}} )
\end{equation} 
be an evolution equation of order $k$ describing $\eta$-pseudo-spherical surfaces. If there exists a local isometric immersion of a surface determined by a solution $u$ for which  the  coefficients of the second fundamental form depend on a jet of finite order of $u$, i.e., $a, b$ and $c$ depend on $x, t, u, \dots, u_{x^{l}}$, where $l$ is finite, then $a, b$ and $c$ are universal, that is $l=0$ and $a, b$ and $c$ depend at most on $x$ and $t$ only.  \end{theorem}

In Section \ref{sec:2}, we give a proof of Theorem \ref{LemCoeff} based on a careful order-by-order analysis of the Codazzi equations (\ref{Eq1}), (\ref{Eq2}) and the Gauss equation (\ref{Gauss}). In Section \ref{sec:3}, we show by means of an example that the class of evolution equations of order $k\geq 3$ for which the components $a,b,c$ are universal in the sense of Theorem \ref{LemCoeff}, that is independent of $u$ and its derivatives, is non-empty.
\section{Proof of the main result}
\label{sec:2}
% Always give a unique label
% and use \ref{<label>} for cross-references
% and \cite{<label>} for bibliographic references
% use \sectionmark{}
% to alter or adjust the section heading in the running head
In the case of a differential equation describing $\eta$-pseudo-spherical surfaces, the structure equations (\ref{struct}) are equivalent to 
\begin{eqnarray} \label{SEq1}
D_t f_{11} - D_x f_{12} = \Delta_{23}\\\label{SEq2}
D_x f_{22} = \Delta_{13}\\\label{SEq3}
D_t f_{31} - D_x f_{32} = -\Delta_{12}
\end{eqnarray}
where $D_t$ and $D_x$ are the total derivative operators and 
\begin{equation}
\Delta_{12} := f_{11}f_{22} -\eta f_{12}; \quad  \Delta_{13} := f_{11}f_{32} - f_{31}f_{12}; \quad  \Delta_{23} = \eta f_{32} - f_{31}f_{22}.
\end{equation}
We shall use the notation 
\begin{equation}
z_{i}=u_{x^{i}}=\frac{\partial^{i}u}{\partial x^{i}},\quad 0\leq i \leq k,
\end{equation}
introduced in \cite{CT} to denote the derivatives of $u$ with respect to $x$ and write the evolution equation (\ref{evolk}) as 
\begin{equation}\label{evolkz}
 z_{0,t} = F(z_0, z_1, \dots, z_k).
\end{equation}
We will thus think of $(t,x,z_{0},\dots,z_{k})$ as local coordinates on an open set of the sub-manifold of the jet space $J^{k}({\bf R}^{2},{\bf R})$ defined by the differential equation (\ref{evolk}).  
We first recall the following lemma from \cite{CT}:
\begin{lemma} Let (\ref{evolkz}) be a $k$-th order evolution equation describing $\eta$-pseudo-spherical surfaces, with associated $1$-forms (\ref{forms}) such that $f_{21}=\eta$. Then necessary conditions for the structure equations (\ref{struct}) to hold are given by
\begin{eqnarray}
f_{11,z_k} = \cdots = f_{11,z_0} = 0\\
f_{21} = \eta\\
f_{31,z_k} = \cdots = f_{31,z_0} = 0\\
f_{12,z_{k}} = 0\\
f_{22, z_k} = f_{22,z_{k-1}} = 0\\
f_{32,z_{k}} = 0\\
f_{11,z_0}^2 + f_{31,z_0}^2 \neq 0
\end{eqnarray}
\end{lemma}
We now proceed with the proof of Theorem \ref{LemCoeff}. 
If $a, b, c$ depend of a jet of finite order, that is  $a, b, c$ are functions of $x, t, z_0, \dots, z_l$ for some finite $l$, then (\ref{Eq1}) and (\ref{Eq2}) become
\begin{eqnarray*}
f_{11}a_t + \eta b_t - f_{12}a_x - f_{22} b_x - 2b \Delta_{13} + (a-c)\Delta_{23}  - \sum_{i=0}^l (f_{12}a_{z_i} + f_{22}b_{z_i})z_{i+1} \\
+\sum_{i=0}^l (f_{11}a_{z_i} + \eta b_{z_i})z_{i,t} = 0,
\end{eqnarray*}
and
\begin{eqnarray*}
f_{11}b_t + \eta c_t - f_{12}b_x - f_{22} c_x + (a-c) \Delta_{13} + 2b\Delta_{23}  - \sum_{i=0}^l (f_{12}b_{z_i} + f_{22}c_{z_i})z_{i+1} \\
+\sum_{i=0}^l  (f_{11}b_{z_i} + \eta c_{z_i})z_{i,t} = 0.
\end{eqnarray*}
Differentiating (\ref{Eq1}) and (\ref{Eq2}) with respect to $z_{l + k}$, and using the fact that $F_{z_k}\neq 0$ and $\eta \neq 0$, it follows that 
\begin{eqnarray}\label{bzlczl}
b_{z_l} = -\frac{f_{11}}{\eta}a_{z_l}, \quad \quad c_{z_l} = \bigg(\frac{f_{11}}{\eta}\bigg)^2a_{z_l}. \
\end{eqnarray}
Differentiating  the Gauss equation (\ref{Gauss}) with respect to $z_l$ leads to  $ca_{z_l} + ac_{z_{l}} - 2bb_{z_{l}} = 0$, and substituting (\ref{bzlczl}) in the latter  leads  to 
\begin{equation}\label{Cases}
\bigg[c+ \bigg(\frac{f_{11}}{\eta}\bigg)^2a + 2\frac{f_{11}}{\eta}b \bigg]a_{z_{l}}=0.
\end{equation}
We therefore have two cases to deal with. The first case corresponds to 
\begin{equation}\label{Case1}
c+ \bigg(\frac{f_{11}}{\eta}\bigg)^2a + 2\frac{f_{11}}{\eta}b \neq  0.
\end{equation} 
It follows then by (\ref{Cases}) that $a_{z_l} = 0$, and hence by (\ref{bzlczl}) that $b_{z_l} = c_{z_l}=0$. Successive differentiating leads to $a_{z_i} = b_{z_i} = c_{z_i} = 0$ for all $i=0, \dots, l$. 
Finally, if the functions $a, b$ and $c$ depend on a jet of finite order, then there are universal, i.e., they are functions of $x$ and $t$ only. We now turn to the second case, defined by the condition
\begin{equation}\label{Case2}
c+ \bigg(\frac{f_{11}}{\eta}\bigg)^2a + 2\frac{f_{11}}{\eta}b = 0,
\end{equation}
on an open set, for which the analysis is far more elaborate. Substituting the expression of $c$ in the Gauss equation $-ac+b^2 = 1$ leads to $ (f_{11}a/ \eta + b )^2 = 1$ so that 
\begin{equation}\label{bandc}
b = \pm 1 - \frac{f_{11}}{\eta}a, \quad  \quad c = \bigg(\frac{f_{11}}{\eta}\bigg)^2 a \mp 2\frac{f_{11}}{\eta}.
\end{equation}
We have then
\begin{eqnarray*}
&D_t b  = -\frac{f_{11}}{\eta}D_t a  - \frac{a}{\eta}f_{11,z_0} F, \quad &D_t c = \bigg(\frac{f_{11}}{\eta}\bigg)^2 D_t a + \frac{2}{\eta}\bigg( \frac{f_{11}}{\eta}a \mp 1\bigg) f_{11,z_0}F,\\
& D_x b  = -\frac{f_{11}}{\eta}D_x a  - \frac{a}{\eta}f_{11,z_0} z_1, \quad 
 & D_x c = \bigg(\frac{f_{11}}{\eta}\bigg)^2 D_x a + \frac{2}{\eta}\bigg( \frac{f_{11}}{\eta}a \mp 1\bigg) f_{11,z_0}z_1,
\end{eqnarray*}
and hence
\begin{eqnarray}
f_{11}D_t a + \eta D_t b & = & -af_{11,z_0}F,\label{1em4}\\
f_{11}D_t b + \eta D_t c & = & \bigg(\frac{f_{11}}{\eta} a \mp 2\bigg)f_{11,z_0}F,\label{2em4}\\
f_{12}D_x a + f_{22}D_x b & = & -\frac{\Delta_{12}}{\eta}D_xa - \frac{af_{22}}{\eta}f_{11,z_0}z_1,\label{3em4}\\
f_{12}D_x b + f_{22}D_x c & = & \frac{f_{11}}{\eta}\frac{\Delta_{12}}{\eta}D_x a  + \frac{\Delta_{12}}{\eta^2} a f_{11,z_0}z_1 + \frac{f_{22}}{\eta}\bigg(\frac{f_{11}}{\eta} a \mp 2\bigg)f_{11,z_0}z_1.
\label{4em4}\end{eqnarray}
Substituting the latter four equalities in (\ref{Eq1}) and (\ref{Eq2}) leads to 
\begin{equation}
-af_{11,z_0}F  +\frac{\Delta_{12}}{\eta}D_xa + \frac{af_{22}}{\eta}f_{11,z_0}z_1- 2b \Delta_{13} + (a-c)\Delta_{23}= 0
\end{equation}
and
\begin{eqnarray*}
\bigg(\frac{f_{11}}{\eta} a \mp 2\bigg)f_{11,z_0}F -  \frac{f_{11}\Delta_{12}}{\eta^{2}}D_x a  -  \frac{\Delta_{12}}{\eta^2} a f_{11,z_0}z_1-\\ \frac{f_{22}}{\eta}\bigg(\frac{f_{11}}{\eta} a \mp 2\bigg)f_{11,z_0}z_1 + (a-c)\Delta_{13} + 2b\Delta_{23} = 0
\end{eqnarray*}
which are equivalent to 
\begin{equation}\label{eq1-1}
-af_{11,z_0}F +\frac{\Delta_{12}}{\eta}a_{x} +\frac{\Delta_{12}}{\eta}\sum_{i=0}^l  a_{z_i}z_{i+1} + \frac{af_{22}}{\eta}f_{11,z_0}z_1- 2b \Delta_{13} + (a-c)\Delta_{23} = 0\\
\end{equation}
and 
\begin{eqnarray}
\bigg(\frac{f_{11}}{\eta} a \mp 2\bigg)\bigg(F  - \frac{f_{22}}{\eta}z_1\bigg)f_{11,z_0}-\frac{f_{11}\Delta_{12}}{\eta^{2}}a_{x} - \label{eq2-1} \\
\frac{f_{11}\Delta_{12}}{\eta^2}\sum_{i=0}^l a_{z_i}z_{i+1}  -  \frac{\Delta_{12}}{\eta^2} a f_{11,z_0}z_1+ (a-c)\Delta_{13}+ 2b\Delta_{23} = 0.
\nonumber
\end{eqnarray} 
We are now led to several cases depending on the value of $l$. 
\begin{itemize}
\item If $l \geq k$, then differentiating (\ref{eq1-1}) with respect to $z_{l + 1}$ leads to $\Delta_{12}a_{z_l}/\eta= 0$. Thus  $a_{z_{l}} = 0$ and also $b_{z_l} = c_{z_l} = 0$ for $l \geq k$ since $\Delta_{12} \neq 0$. 
\item If $l = k-1$, then differentiating (\ref{eq1-1})  and  (\ref{eq2-1}) with respect to $z_k$ leads to 
\begin{eqnarray}\label{eq1-3}
-af_{11,z_0}F_{z_k} + \frac{\Delta_{12}}{\eta} a_{z_{k-1}} = 0,\\\label{eq2-3}
\bigg(\frac{f_{11}}{\eta} a \mp 2\bigg)f_{11,z_0}F_{z_k} - \frac{f_{11}}{\eta}\frac{\Delta_{12}}{\eta}a_{z_{k-1}} = 0.
\end{eqnarray}
Taking into account (\ref{eq1-3}), equation (\ref{eq2-3}) becomes $\mp2f_{11,z_0}F_{z_k} = 0$. 
%Therefore, the system runs into a contradiction when $f_{11,z_0} \neq 0$, which is the case for evolution equations with $f_{ij}$'s satisfying either $H=L=0$, or $H=0$ and $L=0$. If $f_{11,z_0} = 0$, which is the case for evolution equations with $f_{ij}$'s satisfying $H\neq 0$, then $a_{z_{k-1}} = 0$, and hence $b_{z_{k-1}} = c_{z_{k-1}} = 0$
It follows then from (\ref{eq1-3}) that $a_{z_{k-1}} = 0$, and therefore that $b_{z_{k-1}} = c_{z_{k-1}}$ = 0.
\item If $l \leq k-2$, then differentiating (\ref{eq1-1})  and  (\ref{eq2-1}) with respect to $z_k$ leads to 
\begin{eqnarray}\label{eq1-4}
-af_{11,z_0}F_{z_k}  = 0,\\\label{eq2-4}
\bigg(\frac{f_{11}}{\eta} a \mp 2\bigg)f_{11,z_0}F_{z_k}  = 0,
\end{eqnarray}
which imply that 
\begin{equation}\label{f11}
f_{11}=\mu,
\end{equation}
for some real constant $\mu$.
%Again, taking into account (\ref{eq1-3}), equation (\ref{eq2-3}) becomes $\mp2f_{11,z_0}F_{z_k} = 0$. Therefore, the system runs into a contradiction when $f_{11,z_0} \neq 0$. If $f_{11,z_0} = 0$, then $f_{11} = \mu$ and 
Equations (\ref{eq1-1}) and (\ref{eq2-1}) then become
\begin{eqnarray}\label{eq1-5}
\frac{\Delta_{12}}{\eta} a_x+\frac{(\mu f_{22} - \eta f_{12})}{\eta}\sum_{i=0}^l a_{z_i}z_{i+1} - 2b (\mu f_{32} - f_{31}f_{12})\\+ (a-c)(\eta f_{32} - f_{31}f_{22}) = 0\nonumber 
\end{eqnarray}
and
\begin{eqnarray}\label{eq1-55}
-\frac{f_{11}\Delta_{12}}{\eta^{2}}a_x - \frac{\mu(\mu f_{22} - \eta f_{12})}{\eta^2}\sum_{i=0}^l a_{z_i}z_{i+1}  + (a-c)(\mu f_{32} - f_{31}f_{12})\\+ 2b(\eta f_{32} - f_{31}f_{22}) = 0,\nonumber  
\end{eqnarray} 
where
\begin{equation}
\Delta_{12}=\mu f_{22}-\eta f_{12}.
\end{equation}
Note that when $f_{11,z_0} = 0$, the structure equation (\ref{SEq1}) becomes $D_{x}f_{12} = -\Delta_{23}$, or equivalently
\begin{equation}\label{f12f11eq0}
f_{12,z_{k-1}}z_k + \cdots  + f_{12,z_0}z_1 = f_{31}f_{22} - \eta f_{32}. 
\end{equation}
Differentiating (\ref{f12f11eq0}) with respect to $z_k$ leads then to $f_{12,z_{k-1}} = 0$. If  $l = k-2$, then taking into account the latter, and differentiating (\ref{eq1-5}) and (\ref{eq1-55}) with respect to $z_{k-1}$ lead to 
\begin{eqnarray}\label{eq1-6}
\frac{\mu f_{22} - \eta f_{12}}{\eta} a_{z_{k-2}} - 2b \mu f_{32,z_{k-1}} + (a-c)\eta f_{32,z_{k-1}} = 0,\\\label{eq2-6}
- \frac{\mu(\mu f_{22} - \eta f_{12})}{\eta^2} a_{z_{k-2}} + (a-c)\mu f_{32,z_{k-1}} + 2b \eta f_{32,z_{k-1}} = 0.
\end{eqnarray}
Note that $f_{11} = \mu \neq 0$, otherwise we would have $b=\pm 1$ and therefore (\ref{eq2-6}) would become $\pm 2\eta f_{32,z_{k-1}} = 0$ which would lead to a contradiction. Indeed, differentiating the structure equation (\ref{SEq3}) with respect to $z_k$, we obtain $f_{31,z_0}F_{z_{k}} = f_{32,z_{k-1}}$. The vanishing of $f_{32,z_{k-1}}$ would then imply the vanishing of $f_{31,z_0}$, but this is not possible because $f_{11,z_0}^2 + f_{31,z_0}^2 \neq 0$. Therefore, $f_{11} = \mu \neq 0$ and $f_{32,z_{k-1}} \neq 0$. Now, multiplying (\ref{eq1-6}) by $\mu/ \eta$ and adding (\ref{eq2-6}), we now obtain that 
\begin{equation}
\mu (a-c) + 2\eta b =2b\mu^2/\eta   - (a-c) \mu , 
\end{equation}
which is equivalent to 
\begin{equation}\label{acb}
 \eta \mu (a-c) = (\mu^2 - \eta^2) b . 
\end{equation}
If $\mu^2 = \eta^2$, then (\ref{acb}) leads to $a-c =0$ which runs into a contradiction because it follows from (\ref{Case2}), (\ref{f11}) and the hypothesis $\mu^2 = \eta^2$ that $a-c =\mp 2 \mu/\eta \neq 0$.  We have then $\mu^2 - \eta^2 \neq 0$. Substituting (\ref{bandc}) in (\ref{acb}) leads to 
\begin{equation}
\eta \mu \bigg[\bigg(1- \frac{\mu^2}{\eta^2}\bigg)a \pm 2\frac{\mu}{\eta}  \bigg] = (\mu^2 - \eta^2) \bigg[\pm 1 - \frac{\mu}{\eta}a \bigg], 
\end{equation}
which simplifies to $\mu^2 + \eta^2 = 0$ which runs into a contradiction. Finally, if $l < k-2$, where $k\geq 3$, then differentiating (\ref{eq1-5}) and (\ref{eq1-55}) with respect to $z_{k-1}$ and using the non-vanishing of $f_{32,z_{k-1}}$ leads to
\begin{eqnarray}
\eta(a-c) - 2\mu b = 0\\
\mu (a-c) + 2\eta b = 0
\end{eqnarray}
Since $\eta^2 + \mu^2 \neq 0$, we have $a=c$ and $b=0$ which runs into a contradiction with the Gauss equation. 
\end{itemize}
Therefore, for all $l$,  (\ref{Eq1}), (\ref{Eq2}) and the Gauss equation form an inconsistent system. Hence, if the immersion exists then (\ref{Case1}) holds 
and $a$, $b$ and $c$ are functions of $x$ and $t $ only.  
This completes the proof of our theorem. 

\section{An Example}\label{sec:3}

We now show by displaying an example that the class of evolution equations of order $k\geq 3$ for which the components $a,b,c$ are universal in the sense of Theorem \ref{LemCoeff}, that is independent of $u$ and its derivatives, is non-empty. Consider the following fourth-order evolution equation obtained in \cite{FT}
\begin{equation}\label{example}
u_{t}=u_{xxxx}+m_{1}u_{xxx}+m_{2}u_{xx}-uu_{x}+m_{0}u^{2},
\end{equation}  
where $m_{0},m_{1},m_{2}$ are arbitrary real constants. Letting
\begin{equation}
\phi=(m_{1}+2m_{0})u_{xx}+Bu_{x}-\frac{u^{2}}{2}+2m_{0}B,\quad r_{0}=-4m_{0}^{2}B,
\end{equation}
where
\begin{equation}
B=4m_{0}^{2}+2m_{0}m_{1}+m_{2},
\end{equation}
it is straightforward to check that the $1$-forms
\begin{eqnarray}
\omega^1 &=& udx+(u_{xxx}+\phi)dt,\\
\omega^2 &=& -2m_{0}dx+r_{0}dt,\\
\omega^3 &=& udx+(u_{xxx}+\phi)dt,
 \end{eqnarray}
satisfy the structure equations (\ref{struct}) whenever $u$ is a solution of (\ref{example}). Let now 
\begin{equation}
h=e^{2(-2m_{0}x+r_{0}t)},
\end{equation}
and let $\gamma$ and $l$ be real constants such that $l>0$ and $l^{2}>4\gamma^{2}$. The functions $a,b,c$ defined by 
\begin{equation}
a=\sqrt{lh-\gamma^{2}h^{2}-1},\, b=\gamma h,\, c=\frac{\gamma^{2}h^{2}-1}{a},
\end{equation}
satisfy the Gauss equation (\ref{Gauss}) and the Codazzi equations (\ref{Eq1}), (\ref{Eq2}) whenever $u$ is a solution of (\ref{example}).

\begin{acknowledgement}
Research partially supported by NSERC Grant RGPIN 105490-2011 and by the Minist\'erio de Ci\^encia e Tecnologia, Brazil, CNPq Proc. No. 303774/2009-6.
\end{acknowledgement}

\end{document}